\documentclass[12pt]{article}
\usepackage{amsmath,amsthm,amsfonts,amscd,amssymb,pb-diagram}

\theoremstyle{plain}
\newtheorem{Thm}{Theorem}[section]
\newtheorem{Lem}[Thm]{Lemma}
\newtheorem{SDPThm}[Thm]{Semidirect product theorem}
\newtheorem{Crlr}[Thm]{Corollary}
\newtheorem{Prop}[Thm]{Proposition}
\newtheorem*{Rem}{Remark}
\newtheorem*{Conj}{Conjecture}

\newtheorem*{Claim}{Claim}
\theoremstyle{definition}
\newtheorem{Def}[Thm]{Definition}

\def\integ{\mathbb{Z}}

\def\F{\mathcal{F}}

\def\<{\left<}
\def\>{\right>}

\author{Daniel Reuben Krashen}

\title{Severi-Brauer varieties of semidirect product algebras}

\begin{document}

\maketitle

\abstract

A conjecture of Amitsur states that two Severi-Brauer varieties 
are birationally isomorphic if and only if the underlying algebras
are the same degree and generate the same cyclic subgroup of the Brauer
group. It is known that generating the same cyclic subgroup is a necessary
condition, however it has not yet been shown to be sufficient.

In this paper we examine the case where the algebras have a
maximal subfield $K/F$ of degree $n$
with Galois closure $E/F$ whose Galois group is of the form
$C_n \rtimes H$, where $E^H = K$ and $|H|$ is prime to $n$. 
For such algebras we show that the
conjecture is true for certain cases of $n$ and $H$. In particular we prove
the conjecture in the case that $G$ is a dihedral group of order $2p$,
where $p$ is prime.

\section{Introduction}

Let $F$ be a field. We fix for the entire paper a positive integer $n$,
and we suppose that either $n$ is prime, or that $F$ contains a primitive
$n$'th root of unity.
For a field extension $L/F$, and 
$A$ a central simple $L$-algebra, we write $V(A)$ or $V(A/L)$ 
to denote the Severi-Brauer variety of $A$, consisting
of $(deg \ A)$-dimensional right ideals of $A$, and denote the function field
of this variety by $L(A)$. 

We recall the following conjecture:

\begin{Conj}[Amitsur, 1955 \cite{Am}] \label{Amitsur}
Given $A, B$ Central Simple algebras over $F$, $F(A) \cong
F(B)$ iff $[A]$ and $[B]$ generate the same cyclic
subgroup of the $Br(F)$.
\end{Conj}

Amitsur showed in \cite{Am} 
that one of these implications hold. Namely if
$F(A) \cong F(B)$, then the equivalence classes of $A$
and $B$ generate the same cyclic subgroup of the Brauer Group.
The aim of this note is to prove the reverse implication for
certain algebras $A$ and $B$. We will say that the conjecture
holds for the pair $(A, l)$, or simply that $(A, l)$ is true to
mean that $l$ is prime to $deg(A)$ and $F(A) \cong
F(A^l)$. We say that the conjecture is true for $A$ if,
for all $l$ prime to $deg(A)$, $(A, l)$ is true.

One important case is when the algebra $A$ has
a cyclic Galois maximal splitting field. In this case we know that the
conjecture is true for $A$ (\cite{Am}, \cite{Roq}). In this paper we
extend this result to certain $G - H$ crossed products. We recall the
following definitions:
\begin{Def}
Let $G$ be a finite group, and $H$ a subgroup of $G$. A field extension
$K / F$ is called $G - H$ Galois if there exists a field $E$ containing
$K$ such that $E / F$ is $G$-Galois, and $E^H = K$.
\end{Def}
\begin{Def}
Let $A$ be a central simple $F$-algebra. 
$A$ is called a $G - H$ crossed
product if $A$ has a maximal subfield $K$ which is $G - H$ Galois. In the
case $H = 1$, we call $A$ a $G$-crossed product.
\end{Def}
The main theorem in this note
concernes the case of an algebra which
is a so-called semidirect product algebra in the sense of \cite{RoSa:SP}
\begin{Def} \label{SDP def}
$A$ is called a semidirect product algebra if it is a
$G - H$ crossed product where $G = N \rtimes H$.
\end{Def}
This can be interpreted as meaning that $A$ becomes an $N$-
crossed product after extending scalars by some field $K/F$
which is $H$-Galois. In the case where $N$ is a cyclic group, we
will try to exploit the fact that we know Amitsur's conjecture to be
true for $N$ crossed products to prove the conjecture for $G$
crossed products.

\begin{SDPThm} \label{main theorem}
Let $A$ be a semidirect product algebra of degree $n$ as in \ref{SDP def}, 
with $N = C_m = \<\tau\>$, $H = C_n = \<\sigma\>$, such that
the homomorphism $N \rightarrow Aut(H)$ (induced by
conjugation) is injective and $|N|$ and $|H|$ are relatively prime.
Choose $r$ so that we may write $\tau \sigma \tau ^{-1} = \sigma ^r$.
Let
$$S = \frac{\integ[\rho]}{1 + \rho + \rho^2 + \cdots + \rho^{n-1}},$$
and define an action of $\tau$ on $S$ via $\tau(\rho) = \rho^r$, and a ring 
homomorphism $\overline{\epsilon}: S \rightarrow \integ / n \integ$ via 
$\overline{\epsilon}(\rho) = \overline{1}$.
Then $(A,l)$ is true for all $l$ such that 
$\overline{l} \in \overline{\epsilon}\big((S^*)^{\tau}\big)$.
\end{SDPThm}
We give the proof of this in section \ref{proof section}. 
For now, we give the
following corollary:
\begin{Crlr}
Suppose $deg(A) = n$, $n$ odd. If $A$ has a dihedral splitting field 
of degree $2n$. Then the conjecture is true for $A$.
\end{Crlr}
\begin{proof}
In this case, we have $\overline{\epsilon}((S^*)^{\tau})
= (\integ / n \integ)^*$. If $l = 2^k$, $\overline{l}$ is the image of
$(\rho + \rho^{-1})^k$. If $l = 2k + 1$, then $\overline{l}$ is the image
of $\rho^{-k} + \rho^{-(k-1)} + \ldots + \rho^{-1} + 1 + \rho +
\ldots + \rho^{k-1} + \rho^k$. Any other unit in $(\integ / n \integ)^*$
is easily seen to be the image of a product of those above.
\end{proof}

\begin{Rem}
This theorem is already known when $F$ contains the
$n$'th roots of unity, since by a theorem of Rowen and Saltman 
\cite{RoSa:DA}, any such algebra is in fact cyclic, and so the
theorem follows from \cite{Am} or \cite{Roq}.
\end{Rem}

It is worth noting that the hypothesis concerning the splitting field $E$
can be stated in weaker terms for the case $n = p$ a prime number. 
In particular we have:
\begin{Prop}
Suppose $A$ is a central simple $F$-algebra of degree $p$ with a 
maximal subfield $K$, and suppose that there is some extension $E'$ of
$K$ such that $E'/F$ is Galois with group $G = C_p \rtimes H$ where
$(E')^H = K$. Then there is a subfield $E \subset E'$ containing $K$
such that $E/F$ is Galois with group $C_p \rtimes C_m$ where $C_m$ acts
faithfully on $C_p$.
\end{Prop}
\begin{proof}
We define a homomorphism $\phi: H \rightarrow Aut(C_p)$ via the natural
conjugation action of $H$ on $C_p$. Since $Aut(C_p)$ is a cyclic group
every subgroup is cyclic, and we may regard $\phi$ as a surjective
map $H \rightarrow C_m$. Now we define a map
\begin{align*}
G = C_p \rtimes H &\rightarrow C_p \rtimes C_m \\
(a, h) &\mapsto (a, \phi(h))
\end{align*}
one may check quickly that this is a homomorphism of groups and its
kernel is precisely the kernel of $\phi$. Set $H' = ker \phi$, and let
$E = (E')^{H'}$. Since $H'$ is normal in $G$ (as the kernel of a 
homomorphism), we know that $E/F$ is Galois and its Galois group is
$G / H' = C_p \rtimes C_m$. By construction, the action of $C_m$ is faithful
on $C_p$.
\end{proof}
Note also that in the case $n = p$ a prime, $S$ is a ring of 
cyclotomic integers.

\section{Preliminaries}

To begin, let us fix some notation. 
Let $F$ be an infinite field. 
The symbol $\otimes$ when unadorned will always denote
a tensor product over $F$ and $\times$ will denote a fiber product
of schemes over $Spec(F)$. For us an $F$-variety will mean a
quasi-projective geometrically integral seperated scheme of 
finite type over $F$ 
(note $F$ is not assumed to be algebraically closed). By geometrically
integral we mean that the scheme remains integral when fibered up to
the algebraic closure of its field of definition. If $X$ is a variety, we
denote its function field by $F(X)$. We remark that $X$ being
geometrically integral variety implies that $F(X)$ is a regular field
extension of $F$, that is to say, $F(X) \otimes_F F^{alg}$ is a field.

Let $E/F$ be $G$-Galois for some group $G$. If $B$ is an $E$ algebra, 
then a homomorphism $\alpha : G \rightarrow Aut_F(B)$ defines an action
of $G$ on $B$ (as an $F$ algebra) which is called semilinear in case
$$\forall x \in E, b \in B, \ \ 
\alpha(\sigma)(xb) = \sigma(x) \alpha(\sigma)(b).$$
Similarly, if $X$ is an $E$-variety with structure map 
$k: X \rightarrow Spec(X)$, a homomorphism 
$\alpha: G \rightarrow Aut_{Spec(F)}(X)$ defines an action of $G$ on $X$
(as an $F$-scheme) which is called semilinear in case 
$\sigma \circ k = k \circ \alpha(\sigma)$.

Also for $B$ an $E$-algebra as above, given $\sigma \in G$, we define
$^{\sigma} B$ to be the algebra with the same 
underlying set and ring structure as $B$, but with the
structure map $\sigma^{-1} : E \hookrightarrow { ^{\sigma}} B$.
 
Given $A$ a central simple $F$-algebra, we recall that the functor
of points of the Severi-Brauer variety $V(A)$ is given
the following subfunctor of the Grassmannian functor of points
(see \cite{Jah}, \cite{VdB:BS}, or \cite{See:BS}, and \cite{EiHa} for
the definition of the Grassmannian functor):
\begin{equation*}
V(A)(R) = 
\left\{I \subset A_R \left| 
\begin{matrix}
I \text{ is a left ideal and $A_R / I$} \\
\text{is $R$-projective of rank $n$}
\end{matrix}
\right. \right\}
\end{equation*}
and for a homomorphism of commutative $F$-Algebras
$R \overset{\psi}{\rightarrow} S$ we obtain the set map
\begin{align*}
V(A)(\psi) : V(A)(R) &\rightarrow V_k(A)(S) \\
\text{via }I &\mapsto I \otimes_R S
\end{align*}

\subsection{Descent and functors of points}

Given $X$ an $F$-variety, we obtain a functor 
$$X^{E/F}: \{\text{commutative $F$-algebras}\} 
\rightarrow \{\text{sets}\}$$ by
$X^{E/F}(R) = Mor_{sch_E}(Spec(R_E), X_E)$. If $f \in Mor(X,Y)$,
we abuse notation, and refer to the natural transformation induced by
$f$ also by the letter $f$. 

For $\sigma \in G$, we obtain a natural transformation
$\sigma: X^{E/F} \rightarrow X^{E/F}$ via for 
$\phi \in X^{E/F}(R)$,
$\sigma \centerdot \phi = \sigma \circ \phi \circ \sigma^{-1}$.
We denote this action by $\iota_X : G \rightarrow NatAut(X^{E/F})$.

\begin{Prop} \label{functor_descent}
Let $f \in Mor_E(X_E, Y_E)$. Then $f = g_E$ for $g \in Mor_F(X, Y)$ iff
the following diagram commutes:
\begin{equation*}
\begin{diagram}
\node{X^{E/F}} \arrow{e,t}{f} \arrow{s,l}{\sigma}
  \node{Y^{E/F}} \arrow{s,r}{\sigma} \\
\node{X^{E/F}} \arrow{e,t}{f}
  \node{Y^{E/F}}
\end{diagram}
\end{equation*}
Or, in other words, 
For every commutative $F$-algebra $R$ and $\phi \in X^{E/F}(R)$, we have
\begin{equation*}
\sigma \centerdot f(\phi) =f(\sigma \centerdot \phi)
\end{equation*}
\end{Prop}
\begin{proof}
To begin, assume the above condition holds. We have (recalling that 
$f(\phi) = f \circ \phi$),
\begin{gather*}
\sigma \centerdot f(\phi) = \sigma \centerdot (f \circ \phi) =
\sigma \circ (f \circ \phi) \circ \sigma^{-1} \\
f(\sigma \centerdot \phi) = f(\sigma \circ \phi \circ \sigma) =
f \circ \sigma \circ \phi \circ \sigma^{-1}
\end{gather*}
And setting these two to be equal, we have
\begin{equation*}
\sigma \circ f \circ \phi = f \circ \sigma \circ \phi
\end{equation*}
which in turn gives us
\begin{equation*}
f \circ \phi = \sigma^{-1} \circ f \circ \sigma \circ \phi
\end{equation*}
Since this must hold for each $\phi$, this just says that the elements
$f, \sigma^{-1} \circ f \circ \sigma \in Mor_E(X_E, Y_E)$ correspond to the
same natural transformation when thought of as elements of
$Nat(X^{E/F}, Y^{E/F})$ via the Yoneda embedding, 
and therefore they must actually be equal
- that is to say $f = \sigma^{-1} \circ f \circ \sigma$ or
$\sigma \circ f = f \circ \sigma$. But now, by Galois descent of schemes, 
we know $f = g_E$.

Conversely, assume that $f = g_E$. Then by Galois descent of schemes, we 
have $\sigma \circ f = f \circ \sigma$. Now we simply make
our previous argument backwards and find 
that the desired diagram does in fact commute.
\end{proof}

We note the following 
lemma, which can be checked by examining the Grassmannian in terms of
its Pl\"uker embedding:

\begin{Lem}
Suppose that $V$ is an $F$-vector
space, and let $X = Gr_k(V)$. Then the natural semilinear action $\iota_X$
can be described functorially as the natural transformation from
$X^{E/F}$ to itself such that for $R$ an $F$-algebra, $\sigma \in G$,
$M \in X^{E/F}(R)$,
\begin{equation*}
X^{E/F}(\iota_X(\sigma))(M) = \sigma(M) = \{\sigma(m) | m \in M \}
\end{equation*}
where $\sigma$ acts on the elements of $V_{R \otimes E}$ in the natural way.
\end{Lem}

\begin{Crlr} \label{gsbv_galois_action}
Suppose that $A$ is an $F$-central simple algebra, 
and let $X = V_k(A)$. Then the natural semilinear action $\iota_X$
can be described functorially as the natural transformation from
$X^{E/F}$ to itself such that for $R$ an $F$-algebra, $\sigma \in G$,
$I \subset A_{R \otimes E}$ an element of $X^{E/F}(R)$
\begin{equation*}
X^{E/F}(\iota_X(\sigma))(I) = \sigma(I) = \{\sigma(x) | x \in I \}
\end{equation*}
where $\sigma$ acts on the elements of 
$A_{R \otimes E} = A \otimes R \otimes E$ as $id \otimes id \otimes \sigma$.
\end{Crlr}

\subsection{Severi-Brauer varieties of crossed product algebras} \label{xp_sbv}

We give here an explicit birational description of the
Severi-Brauer Variety of a crossed product algebra. An similar discussion
(without the functorial viewpoint) 
may be found in \cite{Sa:LN} (Cor. 13.15). Let $L/F$ be
a $G$-galois extension of degree $n$. 
Let $A = (L, G, c)$ be a crossed product algebra,
where $c$ is taken to be a specific 2-cocycle (not just a cohomology class)
normalized so that $c(id, id) = 1$.

We define the ``functor of splitting 1-chains for $c$'' via:
To begin, define the functor
\begin{gather*}
\F: \left\{
\begin{matrix}
\text{commutative $L$-algebras} \\
\text{with $G$-semilinear action}
\end{matrix}
\right\}
\rightarrow \{\text{sets}\} \\
\F(S) = \{z \in C^1(G, S^*) | \delta z = c \}
\end{gather*}
Where $C^1(G, S^*)$ denotes the set of 1-cochains.
From this we define the functor of splitting 1-chains of $c$ as
\begin{gather*}
Sp_c : \{\text{commutative $F$-algebras}\} \rightarrow \{sets\} \\
Sp_c(R) = \F(R \otimes L)
\end{gather*}

\begin{Prop} \label{crossed_prod_SBVs}
$Sp_c$ is represented by on open subvariety of $U$ of $V(A)$, which is
given as an open subfunctor by
$U(R) = \{I \in V(A)(R) | I + L_R = A_R\}$.
This isomorphism of functors is given by the natural isomorphism
$\Lambda: U \rightarrow Sp_c$, where
$\Lambda(R)(I)$ is the 1-cochain
\begin{gather*}
\sigma \mapsto z(I)_\sigma \\
\text{where } z(I)_\sigma \text{ is the unique element of $L_R$ such that} \\
z(I)_\sigma - u_\sigma \in I.
\end{gather*}
Further, the inverse is given by
\begin{equation*}
\Lambda^{-1}(z) = \underset{\sigma \in G}{\sum} 
(L \otimes R)(z(\sigma) - u_\sigma)
\end{equation*}
\end{Prop}
\begin{proof}
First we note that if $I \in U(R)$ then $I \cap L_R = 0$. This is because
we have the exact sequence of $R$-modules
\begin{equation*}
0 \rightarrow I \cap L_R \rightarrow L_R \rightarrow A_R/I \rightarrow 0
\end{equation*}
which is split since $A_R/I$ is projective. Hence 
$L_R = A_R/I \oplus (I \cap L_R)$, and since $L_R$ is projective (since $L$ is)
we have $I \cap L_R$ is also projective. By additivity of ranks, we have
that $rk(A_R/I) = n$, $rk(L_R) = dim_F(L) = n$ and so $rk(I \cap L_R) = 0$.
Since $I \cap L_R$ is projective, it must be trivial. Consequently,
$I + L_R = A_R$ implies that $A_R = I \oplus L_R$.

To see now that $\Lambda(R)$ is well defined, we just note that
$-u_\sigma \in A_R = I \oplus L_R$, and so there is a unique element
$z(I)_\sigma \in I$ such that $z(I)_\sigma - u_\sigma \in I$. Next we 
check that $z(I)$ defines an element of $Sp_c(R)$. Since $I$ is a left ideal,
\begin{gather*}
z(I)_\tau - u_\tau, z(I)_{\sigma \tau} - u_{\sigma \tau} \in I \\
\Rightarrow 
u_\sigma(z(I)_\tau - u_\tau) - c(\sigma, \tau)(z(I)_{\sigma \tau} - 
u_{\sigma \tau})\in I
\end{gather*}
But this expression may be simplified to:
\begin{multline*}
u_\sigma(z(I)_\tau - u_\tau) - c(\sigma, \tau)(z(I)_{\sigma \tau} - 
u_{\sigma \tau}) + \sigma(z(I)_\tau)(z(I)_\sigma - u_\sigma)\\
= \sigma(z(I)_\tau) u_\sigma -  
c(\sigma, \tau) u_{\sigma \tau} - c(\sigma, \tau) z(I)_{\sigma \tau} + 
c(\sigma, \tau) u_{\sigma \tau}  + \sigma(z(I)_\tau)z(I)_\sigma -
\sigma(z(I)_\tau) u_\sigma \\
= \sigma(z(I)_\tau) z(I)_\sigma - c(\sigma, \tau) z(I)_{\sigma \tau} 
\in I \cap (R \otimes L) = 0
\end{multline*}
and so $\sigma(z(I)_\tau) z(I)_\sigma = c(\sigma, \tau) z(I)_{\sigma \tau}$
which says that 
\begin{equation*}
\delta (z(I))(\sigma, \tau) = \sigma(z(I)_\tau) z(I)_\sigma
(z(I)_{\sigma \tau})^{-1} = c(\sigma, \tau)
\end{equation*}

Next, we check that $\Lambda^{-1}$ is well defined. Let
$z \in Sp_c(R)$ and set $I = \Lambda^{-1}(z)$.
It is clear from the definition that
$I + (L \otimes R) = A_R$, and $A_R / I = 
L \otimes R$ is free (and so projective) of rank $1$. 
To check that $I$ is actually a left ideal, since it follows from the 
definition that $(L \otimes R) I = I$, we need only check
that for each $\sigma \in G$, $u_\sigma I \subset I$, and this in turn
will follow if we can show $u_\sigma (z(\tau) - u_\tau) \in I$ for each 
$\tau \in G$. Calculating, we get
\begin{align*}
u_\sigma \big( z(\tau) - u_\tau \big) 
&= \sigma \big( z(\tau) \big) u_\sigma - 
c(\sigma, \tau)u_{\sigma \tau} \\ 
&= \sigma \big( z(\tau) \big) u_\sigma -
z(\sigma) \sigma \big(z(\tau)\big)z(\sigma \tau)^{-1} u_{\sigma \tau} \\
&= \sigma \big(z(\tau)\big) \big(u_\sigma - 
z(\sigma) z(\sigma \tau)^{-1} u_{\sigma \tau}\big)\\
&= \sigma \big(z(\tau)\big) 
\Big(- \big(z(\sigma) - u_\sigma \big) + \big(z(\sigma) - 
z(\sigma) z(\sigma \tau)^{-1} u_{\sigma \tau}\big) \Big) \\
&= \sigma \big(z(\tau)\big) \Big(- \big(z(\sigma) - u_\sigma \big) 
+ z(\sigma)z(\sigma \tau)^{-1}
\big(z(\sigma \tau) - u_{\sigma \tau} \big) \Big) \\
&= - \sigma \big(z(\tau) \big) \big(z(\sigma) - u_\sigma \big) + 
z(\sigma)\sigma \big(z(\tau) \big) z(\sigma \tau)^{-1} \big(z(\sigma \tau) - 
u_{\sigma \tau} \big) \\
&\in (L \otimes R) \big(z(\tau) - u_\tau \big) + 
(L \otimes R) \big(z(\sigma \tau) -
u_{\sigma \tau} \big) \subset I
\end{align*}
and hence $I$ is a left ideal, and $\Lambda^{-1}$ makes sense.

It remains to show that $\Lambda$ and $\Lambda^{-1}$ are natural 
transformations and are inverses to one another. It follows fairly easily that
if $\Lambda$ is natural and they are inverses of one another then
$\Lambda^{-1}$ will automatically be natural also.

To see that $\Lambda$ is natural, we need to check that for 
$\phi: R \rightarrow S$ a ring homomorphism, and $I \in U(R)$, that
\begin{equation*}
\Lambda(S)(U(\phi)(I)) = Sp_c(\phi)(\Lambda(R)(I))
\end{equation*}
the right hand side is
\begin{equation*}
Sp_c(\phi)(\Lambda(R)(I)) = Sp_c(\phi)(z(I)) = (id_L \otimes \phi)(z(I))
\end{equation*}
and by definition of $z(I)$, we know 
for $\sigma \in G$, $z(I)_\sigma - u_\sigma \in I$
and for the left hand side we have
\begin{equation*}
\Lambda(S)(U(\phi)(I)) = \Lambda(S)(I \otimes_R S)
\end{equation*}
but 
\begin{equation*}
z(I)_\sigma - u_\sigma \in I \Rightarrow 
z(I)_\sigma \otimes 1 - u_\sigma \in I \otimes_R S
\end{equation*}
and now, using the identification
\begin{align*}
(L \otimes R) \otimes_R S &\overset{\sim}{\rightarrow} L \otimes S \\
(l \otimes r) \otimes s &\mapsto l \otimes \phi(r) s
\end{align*}
$z(I)_\sigma \otimes 1$ becomes $(id_L \otimes \phi)(z(I)_\sigma)$, and so
combining these facts gives
\begin{equation*}
Sp_c(\phi)(\Lambda(R)(I))(\sigma) = (id_L \otimes \phi)(z(I)_\sigma) \in
I \otimes_R S \cap (L \otimes S - u_\sigma)
\end{equation*}
and by definition of $\Lambda$, this means 
\begin{equation*}
\Lambda(S)(I \otimes_R S)(\sigma) = (id_L \otimes \phi)(z(I)_\sigma) =
Sp_c(\phi)(\Lambda(R)(I))(\sigma)
\end{equation*}
as desired.

Finally, we need to check that transformations are mutually inverse.
Choosing $I \in U(R)$, we want to show
\begin{equation*}
I = \underset{\sigma \in G}{\sum} (L \otimes R)(z(I)_\sigma - u_\sigma)
\end{equation*}
Now, it is easy to see that the right hand side is contained in the left
hand side. Furthermore, both of these are direct summands of $A_R$ of
corank $n$. For convenience of notation, let us call the right hand side $J$.
\begin{Claim} I / J is projective
\end{Claim}
We show this by considering the exact sequence
\begin{equation*}
0 \rightarrow I / J \rightarrow A_R / J \rightarrow A_R / I \rightarrow 0
\end{equation*}
Since $A_R / I$ is projective, this sequence splits and
$I/J \oplus A_R / I \cong A_R / J$. But since $A_R / J$ is projective,
and $I/J$ is a summand of it, $I/J$ must be projective as well, proving
the claim.

Now, from the exact sequence
\begin{equation*}
0 \rightarrow J \rightarrow I \rightarrow I/J \rightarrow 0
\end{equation*}
we know $rank(I/J) = rank(I) - rank(J) = 0$, and so $I/J = 0$ which says
$I = J$ as desired.

Conversely, if $z \in Sp_c(R)$, we need to verify that
\begin{equation*}
I = \underset{\sigma \in G}{\sum} (L \otimes R)(z(\sigma) - u_\sigma)
\Rightarrow
I \cap (L \otimes R - u_\sigma) = z(\sigma) - u_\sigma
\end{equation*}
But since
$z(\sigma) - u_\sigma \in I$, this immediately follows.
\end{proof}
\begin{Rem}
This same proof will work for an Azumaya algebra (the case where $F$ is a 
commutative ring).
\end{Rem}

This becomes simpler for the case that $L/F$ is a cyclic extension, say
$A = (L/F, \sigma, b)$. In this case, choosing $c$ to be the standard 
2-cocycle:
\begin{equation*}
c(\sigma^i, \sigma^j) =
\begin{cases}
1& i+j < n \\
b& i+j \geq n
\end{cases}
\end{equation*}
If $z \in Sp_c(R)$, then $z$ is determined by its value
on $\sigma$, and 
$z(\sigma)$ must be an element
of $(L \otimes R)$ with ``$\sigma$-norm'' equal to $b$, and conversely it
is easy to check that such an element will determine an element of 
$Sp_c(R)$. With this in mind,
we will write $[N_{L/F} = b]$ for the functor $Sp_c$.
By the above we may write (up to natural isomorphism)
\begin{equation} \label{cyclic_sbv}
[N_{L/F} = b](R) = \{x \in L \otimes R | x \tau(x) \cdots 
\tau ^{m-1}(x) = b \}
\end{equation}
and for a homomorphism $f: R \rightarrow S$, we have:
\begin{equation*}
[N_{L/F} = b](f)(x) = (id_L \otimes f)(x)
\end{equation*}
and by \ref{crossed_prod_SBVs}, this is represented by an open subvariety
of $V(A)$.

\subsection{Group algebra computations} \label{group_algebra}

For convenience of notation, since we will be dealing often with
certain elements of the group algebra $R = \integ G$, we define
for $\gamma \in G$, and $j$ a positive integer
\begin{equation*}
N_{\gamma} ^j = 1 + \gamma + \gamma ^2 + \ldots + \gamma ^{j-1}
\end{equation*}
which we will call the $j$'th partial norm of $\gamma$.

These satisfy the following useful identity which can be easily
verified:
\begin{equation*}
(N_{\gamma} ^j) (N_{\gamma^j} ^i) = N_{\gamma} ^{ij}
\end{equation*}
where $\gamma$ is an element of $G$.

Now, suppose that $u$ is an element in a arbitrary
F-algebra $B$, and $E'$ is a subfield of $B$ such that for all $x
\in E'$, $ux = \gamma(x) u$ for $\gamma \in Aut_F(E')$. Then we
have the identity:
\begin{equation*}
\left (\underset{k=0}{\overset{i-1}{\sum}} 
\left( \underset{j=1}{\overset{k}{\prod}} \gamma^{i-j}(x) \right)
u^{i - k - 1} \right) 
 (x - u)  = \gamma^{i-1}(x) \gamma^{i-2}(x) \cdots
\gamma(x) x - u^i
\end{equation*}
(where we consider the empty product in the case $k=0$ to equal $1$).

If we consider the group algebra $\integ \<\gamma\>$ to act on 
$E'$, then in the above notation, there is an element
$a \in B$ such that
\begin{equation} \label{norm_computation}
a (x - u) = N_{\gamma} ^i x - u^i.
\end{equation}

\subsection{Galois monomial maps}
As in (\ref{cyclic_sbv}), 
let $[N_{E/L} = b]$ be the functor representing elements of norm
$b$.
\begin{Def}
A Galois monomial in $\sigma$ is an element of the group algebra
$\integ\<\sigma\>$.
\end{Def}
Suppose $P$ is a Galois Monomial in $\sigma$. Let
$\epsilon : \integ \left<\sigma\right> \longrightarrow \integ$
be the augmentation map defined by mapping all group elements to
$1$. Then if we set $l = \epsilon(P)$, for every integer $k$, 
$P$ induces a map
\begin{equation*}
P: [N_{E/L} = b^k] \rightarrow [N_{E/L} = b^{kl}]
\end{equation*}
via for $x \in [N_{E/L} = b^k]$,
if $P = \underset{i = 0}{\overset{p-1}{\sum}} n_i \sigma^i$,
\begin{equation*}
P(x) = x^{n_0 + n_1 \sigma + n_2 \sigma^2  + \cdots + n_{p-1} \sigma^{p-1}} = 
\underset{i = 0}{\overset{p-1}{\prod}} \sigma^i(x^{n_i}).
\end{equation*}
We refer to this as the Galois monomial map induced by $P$.

\section{Proof of the semidirect product theorem} \label{proof section}

We begin by fixing notation. Let $A$ be a central simple semidirect product
algebra of degree $n$
as in the statement of theorem \ref{main theorem}, and fix $K/F$
maximal separable in $A$ so that we have the following diagram of fields:
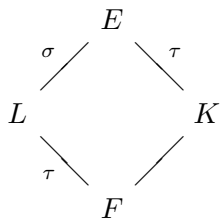
\begin{equation*}
\begin{diagram}
\node[2]{E} \arrow{sw,l,-}{\sigma} \arrow{se,l,-}{\tau}\\
  \node{L} \arrow{se,r,-}{\tau} 
  \node[2]{K} \arrow{sw,-} \\
\node[2]{F}
\end{diagram}
\end{equation*}
Now, as was shown in section \ref{xp_sbv}, 
since $A_L$ is a cyclic algebra,
the functor $[N_{E/L} = b]$ is represented by an open subvariety of 
$V(A_L)$.
The idea of the proving theorem \ref{main theorem}
will be to construct rational maps of Severi-Brauer varieties
by constructing natural transformations between the corresponding
functors.
One way to construct these natural transformations is via Galois monomial
maps.

It can be easily verified that a Galois monomial map
as in the previous section is a natural transformation,
and hence yields an $L$-rational map $V(A_L^k) \rightarrow V(A_L^{kl})$. 
Our goal will be to determine when such a map induces an $F$-rational 
map $V(A^k) \rightarrow V(A^{kl})$. By \ref{functor_descent}, this will
happen when the $\tau$ actions on $[N_{E/L} = b^k] \subset V(A^k)^{E/L}$
and $[N_{E/L} = b^{kl}] \subset V(A^{kl})^{E/L}$
commute with our natural transformation. We will proceed now to determine
the actions of $\tau$, and then to translate these into actions
on the ``norm set'' functors, which will let us answer our question.

\subsection{The Action of $\tau$}

\begin{Lem} \label{action_doesn't_matter}
Let $B = C \otimes L$, where $C$ is a central simple
$F$-algebra of degree $n$. If $\alpha$ is an arbitrary
$\tau$-semilinear action on $B$, then there is an isomorphism
$(B, \alpha) \cong (B, \iota_C)$.
\end{Lem}
\begin{proof}
Let $D=B^\alpha$. Then by descent, 
we have an isomorphism $(D \otimes L, \iota_D) \cong (B, \alpha)$

Since $m=[L : F]$ is relatively prime to $n=deg(B)$, the
restriction map of Brauer groups:
\begin{equation*}
Br_n(F) \overset{res_{L/F}}{\longrightarrow} Br_n(L)
\end{equation*}
is injective. Therefore, since both $D$ and $C$ restrict to the
same element, they are $F$-isomorphic. We can therefore write $D
\cong C$, and again by descent
we get an isomorphism of 
$(D \otimes L, \iota_D) \cong (C \otimes L, \iota_C)$. Combining this
with the isomorphism $(D \otimes L, \iota_D) \cong (B,\alpha)$, 
we have an isomorphism $(C \otimes L, \iota_C) \cong (B, \alpha)$.
\end{proof}
\begin{Crlr}
Let $B$ be a central simple $L$-algebra, $\alpha, \beta$ $\tau$-semilinear
actions on $B$. Then $(B, \alpha) \cong (B, \beta)$.
\end{Crlr}
Therefore, to understand the action of $\tau$ on $A_L$ via $1
\otimes \tau$ up to an isomorphism of pairs, we need only 
define any $\tau$-semilinear action
on $A_L$.

\subsubsection{An Action of $\tau$ on $A_L$}

Since the algebra
$A_L$ has a maximal subfield $E$ which is cyclic over $L$, we may
write $A_L = (E, \sigma, b)$ for some element $b \in L$. Our goal in
this section will be to define a semilinear action of $\tau$ on $A_L$.

Borrowing some of the methods of Rowen and Saltman, we first investigate
the action of $\tau$ on $b \in L$. We first note that since $A$ is an
$F$-algebra, that if we consider the algebra $^\tau A_L$, 
then this is isomorphic to $A_L$ by
the map $id_A \otimes \tau$. On the other hand, one may also check that there
is an isomorphism
\begin{gather*}
^\tau A_L = \ ^\tau (E, \sigma, b) \rightarrow (E, \sigma ^r, \tau(b)) \\
\text{via } E \overset{\tau}{\rightarrow} E \text{ and }
u \rightarrow u
\end{gather*}
and extending to make a homomorphism. Consequently, we have an isomorphism
of central simple algebras $(E, \sigma, b) \cong (E, \sigma^r, \tau(b))$.
In addition, there is also an isomorphism 
$(E, \sigma, b) \cong (E, \sigma^r, b^r)$ (\cite{Pie} p.277 Cor.a), which
means $(E, \sigma^r, b^r) \cong (E, \sigma^r, \tau(b))$. 
This implies $\tau(b) = a b^r$ where 
$a = N_{\sigma^r}(x) = N_\sigma (x)$ for some $x \in E^*$
(\cite{Pie} p.279 Prop.b).

Now to define an action of $\tau$ on $A_L$,
we must first extend the action to the maximal subfield of $A_L$ which
is of the form $L(b^{1/n})$. This will be made more tractable by
choosing a different $b$.
\begin{Lem} 
There exists $b' \in L$
such that $(E, \sigma, b) \cong (E, \sigma, b')$ and such that $\tau(b') = 
\lambda^n (b')^r$ where $\lambda \in L$.
\end{Lem}
\begin{proof}
In the case where $F$ contains the $n$'th roots of unity, this follows
directly from \cite{RoSa:SP}, Lemma 1.2.

For the case where $n = p$ is prime,
We consider the
exact sequence of $\integ / p \integ [\tau]$ modules:
\begin{equation*}
0 \rightarrow \frac{N_{E/L}(E^*)}{(L^*)^p} \rightarrow
\frac{L^*}{(L^*)^p} \overset{\pi}{\rightarrow}
\frac{L^*}{N_{E/L}(E^*)} \rightarrow 0
\end{equation*}
By Maschke's theorem (\cite{Pie}, p.51),
$\integ/p\integ [\tau]$ is a semisimple algebra and
hence every module is projective and every exact sequence splits. We
may therefore choose a splitting map 
$\phi : \frac{L^*}{N_{E/L}(E^*)} \rightarrow \frac{L^*}{(L^*)^p}$.
Let $b'$ be a coset representative for $\phi(b N_{E/L}(E^*))$. 
Since $\phi$ is a splitting, $\pi(b (L^*)^p) = \pi (b' (L^*)^p)$ implies
$\pi(b / b' (L^*)^p) = 1$ which means that $b$ and $b'$ differ by a norm
and so $(E, \sigma, b) \cong (E, \sigma, b')$. Further, since $\phi$
is a $\tau$-morphism, 
\begin{multline*}
\tau(b')(L^*)^p = \tau(\phi(b)(L^*)^p) =
\phi(\tau(b) N_{E/L}(E^*)) = \\ \phi( a b^r N_{E/L}(L^*)) = 
\phi(b^r N_{E/L}(E^*)) = (b')^r (L^*)^p
\end{multline*} 
This gives us $\tau(b') = \lambda^p (b')^r$ for some $\lambda \in L$ as
desired.
\end{proof}
Without loss of generality, we now substitute $b'$ for $b$ and assume
that $\tau(b) = \lambda^n b^r$.

Now consider the field $L(\beta)$, where $\beta$ is defined to be a root
of the polynomial $x^n - b$. We want to show that we can extend the action
of $\tau$ to an order $m$ automorphism of $L(\beta) / F$. To this effect
we first define an map
$\tau' : L(\beta) \rightarrow L(\beta)$, where $\tau'|_L = \tau$ and
$\tau'(\beta) = \lambda \beta^r$.
One may verify this defines an automorphism by considering
$L(\beta) = L[x]/ (x^n - b)$ and noting that $\tau'$ preserves the ideal
$(x^n - b)$. 
\begin{Lem}
We may choose $\lambda$ above so that 
$\tau'$ has order $m$ in $Aut(L(\beta))$.
\end{Lem}
\begin{proof}
Since by definition $\tau' |_L = \tau$, we have 
$(\tau')^m \in Aut(L(\beta)/L)$. We thereby find that 
$ord(\tau') | mn$, $ord(\tau) = m | ord(\tau')$. Therefore we can
write $ord(\tau') = km$, $k | n$, and set $\gamma = (\tau')^k$ and
$M = L(\beta)^{\gamma}$. Since $[M:F] = n$, $[L:F] = m$ have relatively
prime degrees and are both subfields of $L(\beta)$, which has degree $nm$,
we find that $L(\beta) = L \otimes_F M$. Hence we may define
$\tau'' = \tau \otimes id_M \in Aut(L(\beta))$, which is an order $m$
automorphism. But now
$$\tau''|_L = \tau|_L = \tau'|_L$$
and $\tau(b) = \lambda^n b^r \implies \tau''(\beta) = \rho \lambda \beta^r$
where $\rho$ is an $n$'th root of unity. But we see $\tau''$ is defined
in the same way as $\tau'$ except for using $\rho \lambda$ instead of
$\lambda$ as a $n$'th root of unity. Hence, by changing our choice of
$\lambda$ to $\rho \lambda$ we obtain an order $m$ automorphism.
\end{proof}

For simplicity of notation we denote the extension $\tau'$ 
of $\tau$ to $L(\beta)$
also by $\tau$.
By the above description, we have
\begin{equation*}
\tau \beta = \lambda \beta ^r
\end{equation*}
where $\lambda \in L$.

We now use this information to define an action of $\tau$ on $A$.
Since $A_L = (E, \sigma, b)$ can be thought of as the free noncommutative 
$F$-algebra generated by $E$ and $u$ modulo the relations 
$u x  - \sigma (x) u = 0$ and $u^n = b$, giving an $F$-homomorphism 
$A_L \rightarrow B$
is equivalent to giving an $F$-map $\phi: L \rightarrow B$ and choosing
and element $\phi(u) \in B$ such that $\phi(u) \phi(x) - \phi(\sigma(x))
\phi(u) = 0$ and $\phi(u)^n - \phi(b) = 0$. 
Consequently, since any $F$-endomorphism of $A_L$
is an automorphism (since $A_L$ is finite dimensional and simple), to define
an action of $\tau$ on $A_L$, we need only define $\tau$ on $E$ and on
$u$ and then check that our relation is preserved.

To begin, we define $\tau |_E : E \rightarrow E \subset A_L$ 
to be the original Galois action, and $\tau(u) = \lambda u^r$. Checking
our relations we have:
\begin{multline*}
\tau(u) \tau(x) = \lambda u^r \tau(x)
= \lambda \sigma ^r \tau(x) u^r
= \lambda \tau \sigma(x) u^r
= \tau \sigma (x) \lambda u^r \\
=\tau (\sigma(x)) \tau(u)
\end{multline*}
and
\begin{equation*}
\tau(u^n) = \lambda^n u^{rn} = \lambda^n b^r = \tau(b)
\end{equation*}
Since
$L(u) \cong L(\beta)$ where $L(\beta)$ is as above, we know that 
$\tau^m(u) = u$. Since $\tau$ has order $m$ on $E$, together this means
that $\tau$ as defined above is an order $m$ semilinear
automorphism of $A_L$. We will refer to this action as 
$\alpha : \<\tau\> \rightarrow Aut_F(A_L)$.

\subsubsection{An Action of $\tau$ on $A_L ^l$}

We define $A_L ^l \subset \otimes^l A_L =
\overset{\text{$l$ - times}}{\overbrace{A_L \otimes_L \ldots
\otimes_L A_L}}$ to be the algebra generated by $E \otimes_L 1
\otimes_L \ldots \otimes_L 1$ (which we will identify with just
$E$), and $v = u \otimes_L u \otimes_L \ldots \otimes_L u$.
\begin{Lem}
$[A_L ^l] = [A_L]^l$, where brackets denote classes in $Br(L)$
\end{Lem}

\begin{proof}
Since $A_L \cong (E, \sigma, b)$, we simply need to verify that
$A_L ^l$ is just the symbol algebra $(E, \sigma, b^l)$. But this
follows because we clearly have $A_L ^l = \overset{p}{\underset{i
= 0}{\coprod}} Ev^i$, and we need only check the two defining
identities:
\begin{align*}
v^p = (u \otimes u \otimes \ldots \otimes u)^p &= u^p \otimes u^p
\otimes \ldots \otimes u^p \\
&= b \otimes b \otimes \ldots \otimes b \\
&= b^l \otimes 1 \otimes \ldots \otimes 1
\end{align*}
and
\begin{align*}
v x &= (u \otimes u \otimes \ldots \otimes u) (x \otimes 1 \otimes
\ldots \otimes 1) \\
&= (u x \otimes u \otimes \ldots \otimes u) \\
&= (\sigma(x) u \otimes u \otimes \ldots \otimes u) \\
&= (\sigma(x) \otimes 1 \otimes \ldots \otimes 1) (u \otimes u
\otimes \ldots \otimes u) \\
&= \sigma(x) v
\end{align*}
\end{proof}
Next we note that we have a $\tau$-semilinear action on
$\otimes^l A_L$ which is induced (diagonally) by the $\tau$ action
on $A_L$, and further, since it is easy to establish that:
\begin{equation*}
\tau v = \lambda^l v^r
\end{equation*}
and the action of $\tau$ is the usual one on $E$, we know that
$A_L ^l$ is preserved by $\tau$ and hence we have an induced
action on $A_L ^l$. We call this action $\alpha^l$

\subsubsection{$\tau$-Action on Norm Sets}

Our goal now will be to describe an action on the norm sets which is
compatible with the above $\tau$-action on ideals. The following lemma assures
us that since the actions on the algebras $A_L$ and $A_L^l$
given above are isomorphic to the standard actions, they also induce 
isomorphic actions on $V(A)^{L/F}$ and $V(A^l)^{L/F}$ 
respectively. Therefore, we may proceed to find actions on
the norm sets compatible with
the $\tau$ actions given above.

Let $(B, \beta)$ be an algebra with $\tau$-semilinear action
such that $B$ central
simple over $L$. Then by corollary \ref{gsbv_galois_action},
we have an induced action on 
$V(B)(R \otimes L)$ via for $I \in V(B)(R \otimes L)$,
thinking of $I \subset B \otimes_L (L \otimes R)$
\begin{equation*}
\beta (\tau) I = \{ \beta(\tau)(x) | x \in I \}
\end{equation*}
where $\beta(\tau)$ is acting here on 
$B \otimes_L (L \otimes R) = B \otimes R$ as 
$\beta(\tau) \otimes 1$.
\begin{Lem}
If $f: (B, \beta) \rightarrow (B', \beta')$ is an isomorphism, then
the induced isomorphism $V(B)(\_) \rightarrow V(B')(\_)$
commutes with the actions of $\tau$
\end{Lem}
\begin{proof}
This is a simple check:
\begin{multline*}
f(\beta(\tau)(I)) = f( \{ \beta(\tau)(x) | x \in I \}) =
\{ f(\beta(\tau)(x)) | x \in I \} = \\
\{\beta'(\tau)(f(x)) | x \in I \} = \beta'(\tau)(f(I))
\end{multline*}
\end{proof}
As was shown in the previous section, we know $\tau$ acts on $A_L$
and hence also on $A_L \otimes R$. This translates to $\tau$
acting on $E \otimes R$ via the natural action on $E$, 
and by $\tau u = \lambda
u^r$ (where we abuse notation by writing $u$ for $u \otimes 1$). 
Therefore, if $x-u \in I$ then $\tau x - \tau u \in \tau
I$, or in other words, $\tau x / \lambda - u^r \in \tau I$. To translate
this into an action on norm sets, we recall that our birational identification
between $V(A_L)$ and $[N_{E/L} = x]$ is via $I \in V(A_L)^{E/L}(R)$
being identified with $I \cap (E_R - u)$. Therefore to find our $\tau$ action
on norm sets, we take an ideal $I$ with a given intersection $x - u \in
E_R - u$ and find the intersection of the new ideal $\tau(I)$.

Since $GCD(r, n) = 1$ (because $r^m \equiv_n 1$), 
we may select a positive integer $t$ so
that $rt = sn + 1$. By equation \ref{norm_computation} in 
\ref{group_algebra}, there is
an $a \in A$ such that
\begin{align*}
a(\tau x / \lambda - u^r) &= N_{\sigma^r} ^ t(\tau x / \lambda) - u^{rt} \\
&= N_{\sigma^r} ^t (\tau x / \lambda) - b^s u \\
&= \tau N_{\sigma} ^t (x) / N_{\sigma^r} ^t (\lambda) - b^s u \\
&= \tau N_{\sigma} ^t (x) / \lambda ^t - b^s u
\end{align*}
where the last step follows from the fact that $\lambda \in L =
E^\sigma$.

Now, since $\tau I$ is a left ideal containing $\tau x / \lambda
- u^r$, it must also contain $\tau N_{\sigma} ^t (x) / \lambda ^t
b^s - u$. Therefore, $\tau N_{\sigma} ^t (x) / \lambda ^t b^s - u
\in \tau I \cap (E - u)$. This tells us precisely that $\tau I$
corresponds to $\tau N_{\sigma} ^t (x) / \lambda ^t b^s$,
and so we get an action of $\tau$ on
$[N_{E/L} = b]$ via
\begin{equation*}
x \overset{\tau}{\mapsto} \tau N_{\sigma} ^t (x) / \lambda ^t b^s
\end{equation*}
which makes the following diagram commute:
\begin{equation*}
\begin{diagram}
\node{V(A)^{E/L}} \arrow{e,t}{\tau}
  \node{V(A)^{E/L}} \\
\node{[N_{E/L} = b]} \arrow{n} \arrow{e,t}{\tau}
  \node{[N_{E/L} = b]} \arrow{n}
\end{diagram}
\end{equation*}

Similarly, using the fact that $\tau$ acts on $A^l$ via $v
\mapsto \lambda ^l v^r$, $v^p = b^l$, and $v$ induces $\sigma$ on
$E$, we get that $\tau$'s action on $V(A^l)$ yields an action on
$[N_{E/L} = b^l]$ via:
\begin{equation*}
x \overset{\tau}{\mapsto} \tau N_{\sigma} ^t (x) / \lambda ^{lt}
b^{ls}.
\end{equation*}

Now, suppose $P \in \integ\<\sigma\>$, $\epsilon(P) = l$.
As $P$ can be considered as a map from $[N_{E/L} = b^k]$ to
$[N_{E/L} = b^{k l}]$, it is acted upon by $\tau$ via $\tau \bullet P
= \tau \circ P \circ \tau ^{-1}$. On the other hand, there is a
natural action of $\tau$ on the group algebra $\integ\<\sigma\>$
given by conjugation by $\tau$, i.e. $\sigma
\overset{\tau}{\mapsto} \sigma^r$. We claim that these two
actions coincide. Thinking of $P$ as an element of the
group algebra $\integ G$, we write $\tau P \tau^{-1}$ as the
action of $\tau$ by conjugation (as a group element).

\begin{Prop}
$\tau \bullet P = \tau P \tau^{-1}$
\end{Prop}
\begin{proof}
We aim to show
$\tau \circ P = \tau P \tau^{-1} \circ \tau$.
To start choose $x \in [N_{E/L} = b^k](R)$. Since
$P(x) \in [N_{E/L} = b^{k l}](R)$,
we have
\begin{align*}
\tau \circ P(x) &= \tau N_{\sigma} ^t (P(x)) / \lambda ^{klt} b^{kls} \\
&= \tau P N_{\sigma} ^t (x) / \lambda ^{klt} b^{kls}
\end{align*}
On the other hand,
\begin{align*}
\tau P \tau^{-1} \circ \tau (x) &= \tau P \tau^{-1} (\tau
N_{\sigma}
^t (x) / \lambda ^{kt} b^{ks}) \\
&= \tau P N_\sigma ^t (x) / \tau P \tau^{-1} (\lambda^{kt} b^{ks}) \\
&= \tau P N_\sigma ^t (x) / \lambda^{klt} b^{kls} \\
&= \tau \circ P(x)
\end{align*}
Where the second to last step follows from the fact that $\tau P \tau^{-1}$
is a monomial in $\sigma$, and that $\lambda$ and $b$ are
$\sigma$-fixed. To finish, we see that by composing on the right
by the map $\tau$, we get
\begin{equation*}
\tau P \tau^{-1} = \tau \circ P \circ \tau^{-1} = \tau \bullet P
\end{equation*}
as desired.
\end{proof}

\begin{Crlr} \label{fixed_monomials_commute_with_tau}
If $P \in \integ\left<\sigma\right> ^\tau$ then the induced map on
norm sets 
\begin{equation*}
P: [N_{E/L} = b^k] \rightarrow [N_{E/L} = b^{kl}]
\end{equation*}
commutes with the action of $\tau$.
\end{Crlr}

\subsection{Proof of the Main Theorem}

We recall our earlier definitions of $\epsilon$,
$\overline{\epsilon}$, and $S$.

\begin{Lem}
Let $P_1, P_2 \in \integ\<\sigma\>$, $\epsilon(P) = l_i$. Then 
for any $k \in \integ$, $P_1 P_2 $ induces a map
\begin{equation*}
P_1 P_2: [N_{E/L} = b^k] \rightarrow [N_{E/L} = b^{k l_1 l_2}]
\end{equation*}
which is the composition of the maps
\begin{align*}
P_2: [N_{E/L} = b^k] &\rightarrow [N_{E/L} = b^{k l_2}] \\
\text{and }\,\, P_1: [N_{E/L} = b^{k l_2}] &\rightarrow [N_{E/L} =
b^{k l_1 l_2}]
\end{align*}
\end{Lem}
\begin{proof}
This just comes from the fact that the group algebra acts on
$E^*$ with composition being identified with multiplication in
the group algebra.
\end{proof}

\begin{Def}
For $i, k \in \integ$, we define a natural transformation (morphism)
\begin{equation*}
\phi_k : [N_{E/L} = b^i] \rightarrow [N_{E/L} = b^{i + nk}]
\end{equation*}
by the rule: for $x \in [N = b^i](R)$, $\phi_j(x) = xb^k$.
\end{Def}

Note that we abuse notation here, and don't specify the domain or
range of $\phi_j$ in the notation. In any case, one may easily
verify that $\phi_j \circ \phi_k = \phi_{j+k}$. In particular,
these maps are all invertible and hence are birational morphisms.

\begin{Lem}
$\phi_k$ commutes with the action of $\tau$.
\end{Lem}
\begin{proof}
We consider $\phi_k : [N_{E/L} = b^i] \rightarrow [N_{E/L} = b^{i
+ nk}]$. Using the formulas for the $\tau$ actions described
earlier, we have:
\begin{align*}
\tau(\phi_k(x)) &= \tau(xb^k) \\
&=\frac{\tau N_{\sigma} ^t (b^k x)}{\lambda^{(i + nk)t} b^{(i
+nk)s}} \\
&=\frac{\tau N_{\sigma} ^t (b^k) \tau N_{\sigma} ^t
(x)}{\lambda^{it} b^{is} \lambda^{nkt} b^{nks}} \\
&=\frac{\tau (b^{tk}) \tau N_{\sigma} ^t (x)}
{\lambda^{it} b^{is} \lambda^{nkt} b^{nks}} \\
&=\frac{(\lambda^n b^r)^{tk} \tau N_{\sigma} ^t (x)}
{\lambda^{it} b^{is} \lambda^{nkt} b^{nks}} \\
&=\frac{\lambda^{nkt} (b^{sn + 1})^k \tau N_{\sigma} ^t (x)}
{\lambda^{it} b^{is} \lambda^{nkt} b^{nks}} \\
&=\frac{b^k \tau N_{\sigma} ^t (x)}
{\lambda^{it} b^{is}} \\
&=b^k \tau(x) \\
&= \phi_k (\tau(x))
\end{align*}
\end{proof}

\begin{Lem} \label{commutation}
For $P$ any Galois monomial in $\sigma$, 
$P \phi_k = \phi_{\epsilon(P)k} P$
\end{Lem}
\begin{proof}
If $P = \underset{i = 0}{\overset{n-1}{\sum}} n_i \sigma^i$ then we simply
compute:
\begin{multline*}
P \circ \phi_k(x) = P(b^k x) = 
\underset{i = 0}{\overset{n-1}{\prod}} (\sigma^i(b^k x))^{n_i} = 
\underset{i = 0}{\overset{n-1}{\prod}} 
(\sigma^i(b^k))^{n_i} (\sigma^i(x))^{n_i} \\ =
\underset{i = 0}{\overset{n-1}{\prod}} 
(b^k)^{n_i} 
\underset{i = 0}{\overset{n-1}{\prod}} 
(\sigma^i(x))^{n_i} = (b^k)^{\sum n_i} P(x) = b^{\epsilon(P) k} P(x) 
= \phi_{\epsilon(P) k} P(x)
\end{multline*}
\end{proof}
We now prove the main theorem:
\begin{Thm}
If $l \in \integ$ such that $\overline{l} \in
\overline{\epsilon}((S^*)^{\tau})$ then there is a birational map
\begin{equation*}
V(A) \rightarrow V(A^l)
\end{equation*}
\end{Thm}
\begin{proof}
We note first that $\integ\left<\sigma\right>^\tau \rightarrow S^\tau$ is surjective,
since if we consider the short exact sequence:
\begin{equation*}
0 \rightarrow N\integ \rightarrow \integ\<\sigma\> \rightarrow S \rightarrow 0
\end{equation*}
We get a long exact sequence in group cohomology
\begin{equation*}
N\integ \rightarrow \integ\<\sigma\>^\tau \rightarrow S^\tau \rightarrow
H^1(\tau, N\integ)
\end{equation*}
and since as a $\tau$-module, $N\integ \cong \integ$ with trivial action,
we have
\begin{equation*}
H^1(\tau, N\integ) = Hom(\left<\tau\right>, \integ) = 
Hom(\integ / m \integ, \integ) = 0
\end{equation*}
giving us a surjective map $\integ\<\sigma\>^\tau \rightarrow S^\tau$ as
claimed.

Now, choosing $\alpha \in (S^*)^\tau$ with $\overline{\epsilon}(\alpha) = 
\overline{l}$, we can find $\widetilde{\alpha} \in 
\integ\left<\sigma\right>^\tau$ mapping to $\alpha$. If 
$\epsilon(\widetilde{\alpha}) = l'$, then by the commutativity of the diagram
\begin{equation} \label{big_diagram}
\begin{diagram}
\node{0} \arrow{e}
 \node{\integ N} \arrow{e}
  \node{\integ\left<\sigma\right>} \arrow{e} \arrow{s,l}{\epsilon}
   \node{S} \arrow{e} \arrow{s,l}{\overline{\epsilon}}
    \node{0} \\
  \node[3]{\integ} \arrow{e}
   \node{\integ / n \integ}
\end{diagram}
\end{equation}
we have $\overline{l'} = \overline{l}$ and so $l' = l + kn$ for some 
$k \in \integ$. 
This means that $\widetilde{\alpha} : [N_{E/L} = b] \rightarrow
[N_{E/L} = b^{l + nk}]$, which commutes with the action of $\tau$ by
\ref{fixed_monomials_commute_with_tau}. Composing this with the map 
$\phi_{-k}$ gives us a natural transformation
$\phi_{-k} \circ \widetilde{\alpha} : 
[N_{E/L} = b] \rightarrow [N_{E/L} = b^l]$ 
which commutes with $\tau$.
We will show that this actually induces an natural isomorphism.

Next pick
$\beta \in (S^*)^\tau$ such that $\alpha \beta = 1 \in S$, and choose
$\widetilde{\beta} \in \integ \left<\sigma\right>^\tau$ mapping to $\beta$,
and let $\big(\epsilon({\widetilde{\beta}})\big) l = 1 + ns$. 
We now have
$\widetilde{\beta} : [N_{E/L} = b^l] \rightarrow [N_{E/L} = b^{1 + ns}]$ and
composing with $\phi_{-s}$ yields a natural transformation
$ \phi_{-s} \circ \widetilde{\beta} : 
[N_{E/L} = b^l] \rightarrow [N_{E/L} = b]$ 
which commutes with $\tau$.

Now, we compose $\phi_{-k} \circ \widetilde{\alpha}$ with 
$\phi_{-s} \circ \widetilde{\beta}$, which by construction is a natural
transformation $[N_{E/L} = b] \rightarrow [N_{E/L} = b]$. 
If we write $\widetilde{\alpha} \widetilde{\beta} = 1 + rN$, then
using \ref{commutation}, we compute:
\begin{equation*}
\phi_{-s} \circ \widetilde{\beta} \circ \phi_{-k} \circ \widetilde{\alpha}
= \phi_{-s} \phi_{-\epsilon(\widetilde{\beta}) k} \widetilde{\beta}
\widetilde{\alpha} = \phi_{-(s + \epsilon(\widetilde{\beta}) k)}(1 + rN)
\end{equation*}
But $1 + rN : [N_{E/L} = b] \rightarrow [N_{E/L} = b^{1 + rn}]$ 
is simply the map $\phi_r$, so
\begin{equation*}
\phi_{-s} \circ \widetilde{\beta} \circ \phi_{-k} \circ \widetilde{\alpha}
= \phi_{-(s + \epsilon(\widetilde{\beta}) k)} \phi_r
= \phi_{r - s - \epsilon(\widetilde{\beta}) k}
\end{equation*}
which is clearly an isomorphism (one can check in fact that 
$r - s - \epsilon(\widetilde{\beta}) k = 0$ giving that the right hand 
side above is the identity). This argument shows that
$ \phi_{-k} \circ \widetilde{\alpha} : 
[N_{E/L} = b] \rightarrow [N_{E/L} = b^l]$ is an also an isomorphism 
which therefore induces a birational map
\begin{equation*}
V(A) \rightarrow V(A^l)
\end{equation*}
\end{proof}

\bibliographystyle{alpha}
\bibliography{citations}

\end{document}